\documentclass{article}
\usepackage{amsfonts,amsthm,amsmath,amssymb,enumerate,geometry,latexsym,mathrsfs}
\usepackage[all]{xy}
\usepackage{enumerate}

\newcommand{\ud}{\mathrm{d}}

\newcommand{\cat}{\mathrm{cat}}
\newcommand{\ccat}{\mathrm{ccat}}

\newcommand{\intl}{\mathrm{int}}

\newcommand{\im}{\mathrm{Im}}

\newcommand{\mon}{\mathrm{Mon}}
\newcommand{\homo}{\mathrm{Hom}}

\newtheorem{The}{\textbf Theorem}[section]
\newtheorem{Pro}[The]{\textbf Proposition}
\newtheorem{Cor}[The]{\textbf Corollary}
\newtheorem{Lem}[The]{\textbf Lemma}

\newtheorem{Not}[The]{\textbf Notation}

\theoremstyle{definition}
\newtheorem{Def}[The]{Definition}

\theoremstyle{remark}
\newtheorem{Rem}[The]{\bf Remark}

\theoremstyle{definition}
\newtheorem{Exp}[The]{Example}

\def\bproof{\textbf{Proof}: }
\def\eproof{\hfill$\Box$}

\usepackage[dvips]{graphicx}

\date{}
\author{Tieqiang Li \\Department of Mathematical Sciences\\Durham University\\Durham, England \\tieqiang.li@durham.ac.uk\and Dirk Sch\"utz \\Department of Mathematical Sciences\\Durham University\\Durham, England\\dirk.schuetz@durham.ac.uk}

\title{\bf On the relative Lusternik-Schnirelmann category with respect to a closed 1-form}
\begin{document}
\maketitle

\begin{abstract}
In this article we study a homotopy invariant $\cat(X,B,[\omega])$ on a pair $(X,B)$ of finite CW complexes with respect to a continuous closed 1-form $\omega$. This is a generalisation of a Lusternik-Schnirelmann category type $\cat(X,[\omega])$, developed by Farber in \cite{Far1, Far2}, studying the topology of a closed 1-form. The article establishes the connection with the original notion $\cat(X,[\omega])$ and obtains analogous results on critical points and homoclinic cycles. We also provide a similar ``cuplength'' lower bound for $\cat(X,B,[\omega])$.
\end{abstract}

\begin{center}              \section{\sc Introduction}                      \end{center}
Michael Farber \cite{Far1,Far2} initiated a systematic study of a generalisation of the classical Lusternik-Schnirelmann category with respect to a closed 1-form, $\cat(X,\xi)$, on a finite CW complex $X$. In \cite {Far1} the power of such a notion is demonstrated in the study of the topology of critical points and the existence of homoclinic cycles on a closed manifold. Compared to the Morse inequalities of a Morse closed 1-form, $\cat(X,\xi)$ is applicable to more degenerate conditions, but in general it is harder to compute. In \cite{FS1,FS2} Farber and Sch\"utz improve the previous results and give more detailed insights on this issue.

In this article we generalise the controlled version of the above notion to the relative case on a finite CW pair $(X,B)$, which coincides with the absolute one when the subset $B$ is empty. In particular, Section 2 introduces the definition of this relative category $\cat (X,B,\xi)$, and in Section 3 we describe the immediate properties of the object. As a main result, we obtain the inequality relating the relative category with the absolute ones. We summarise this in the following theorem:

\begin{The}
Suppose $X$ is a finite CW-complex and $A,B$ are subcomplexes of $X$ with $A\subset B$, and let $\xi\in H^1(X,\mathbb{R})$ be a cohomology class of $X$ and $i^\ast:H^1(X;\mathbb{R})\to H^1(B;\mathbb{R})$ be the induced map of the inclusion map $i:B\to X$, then we have the following inequality: 
\begin{eqnarray*}
\cat(X,A,\xi) & \leq & \cat(X,B,\xi)+ \cat \big( B,A,i^{\ast}(\xi) \big).
\end{eqnarray*}
\end{The}

Note that $\xi$ needs not restrict to the trivial cohomology class on $B$. In the case of $\xi=0$, $\cat(X,A,\xi)$ reduces to the usual relative Lusternik-Schnirelmann category, and this result is given in \cite{Cornea}.

In Section 4 we relate this relative Lusternik-Schnirelmann category to the existence of homoclinic cycles for gradient-like vector field on a manifold with boundary, generalizing previous work of Farber \cite{Far1}.

\begin{The}
Let $M$ be a smooth compact manifold with boundary $\partial M$, and $\omega$ be a closed 1-form on $M$ satisfying certain transversality conditions on the exit set $B\subset \partial M$. If the number of critical points of $\omega$ is less than $\cat(M,B,[\omega])$, then any gradient of $\omega$ transverse on $(\partial M,B)$ contains at least one homoclinic cycle.
\end{The}

The transversality conditions above prescribe a ``nice'' behaviour near the boundary $\partial M$, which is explained in more detail in Section 4. In particular, the exit set $B$ is a $0$-codimensional submanifold of $\partial M$ possibly with boundary.

\begin{center}                   \section{\sc Definition of $\cat (X,B,\xi)$}                                        \end{center}
Firstly, we recall the definition for closed 1-forms on topological spaces resembling the essential features of the conventional closed 1-forms in differential topology. This is first defined in \cite{Far1}.

\begin{Def}
Let $X$ be a topological space, {\it a continuous closed 1-form $\omega$ on $X$} is defined to be a collection $\{f_U\}_{U\in \mathscr{U}}$ of continuous real functions $f_U: U\to \mathbb{R}$, where $\mathscr{U}=\{U\}$ is an open cover of $X$ such that for any pair $U,V\in \mathscr{U}$, the difference
$$f_U|_{U\cap V}-f_V|_{U\cap V}: U\cap V\to \mathbb{R}$$
is locally constant.
\end{Def}
In Chapter 10.2 of \cite{Far}, Farber provides a comprehensive description of this notion, here we only recollect the essential properties necessary for our study.

Two continuous closed 1-forms $\omega_1=\{f_U\}_{U\in \mathscr{U}},\omega_2=\{g_V\}_{V\in \mathscr{V}}$ are called {\it equivalent} if the union $\{f_U,g_V\}_{U\in \mathscr{U},V\in \mathscr{V}}$ of the collections is a continuous closed 1-form, i.e. for any $U\in \mathscr{U}$ and $V\in \mathscr{V}$, the differnce $f_U-g_V$ of the two functions $f_U,g_V$ is locally constant on $U\cap V$. A trivial example for such topological continuous closed 1-form can be constructed as follows:

\begin{Exp}
Suppose we take the whole space $\{X\}$ as the open cover, then any continuous function $f:X\to \mathbb{R}$ defines a continuous closed 1-form on $X$, denoted as $\ud f$. It can be seen as the continuous version of an exact form in differential topology, and we call it {\it continuous exact 1-form}.

In such an example, two exact 1-forms $\ud f,\ud g$ are equivalent $\ud f=\ud g$ if and only if $f-g:X\to \mathbb{R}$ is locally constant, i.e. constant on each connected component of $X$.
\end{Exp}

\begin{Exp}
Consider the 1-dimensional sphere $S^1$ parametrized by $t\to e^{\pi it}$ and cover it with ${U,V}$ where $U=(-\frac{1}{6},\frac{7}{6})$ and $V=(\frac{5}{6},\frac{13}{6})$. Let $\theta_U$ and $\theta_V$ be angular functions, i.e. $\theta_U(x)=\pi x$ for $x\in U$ and $\theta_V(y)=\pi y$ for $y\in V$. Then $\theta_V|_{U\cap V}-\theta_U|_{U\cap V}$ is locally constant, hence $\ud \theta=\{\theta_U,\theta_V\}$ is a continuous closed 1-form on $S^1$. It is easy to see that $\ud \theta$ is not exact.
\end{Exp}

We want to define integration for topological closed 1-forms, which leads to the cohomology class.
\begin{Def}
Suppose we have a closed 1-form $\omega=\{f_U\}_{U\in \mathscr{U}}$ for some open cover $\mathscr{U}=\{U\}$ of topological space $X$, and $\gamma:[0,1]\to X$ is a continuous path on $X$. The line integral $\int_\gamma \omega$ is defined as follows:
\begin{eqnarray*}
\int_{\gamma}\omega & = & \sum_{i=0}^{n-1} \left( f_{U_i}(\gamma(t_{i+1}))-f_{U_i}(\gamma(t_i))\right),
\end{eqnarray*}
where $t_0=1<t_1<\cdots<t_n=1$ is a partition of the closed inteval $[0,1]$ such that $\gamma[t_i,t_{i+1}]\subset U_i$ with $U_i\in \mathscr{U}$, for all $1\leq i\leq n$.
\end{Def}

\begin{Rem}
This integration is independent of the choice of partitions and the open cover $\mathscr{U}$, and only depends on the homology class of the path relative to its end points, see \cite[\S 10.2]{Far}.
\end{Rem}



\begin{Def}
Let $\omega$ be a closed 1-form on a topological space $X$, the {\it homomorphism of periods}: $\pi_1(X,x_0)\to \mathbb{R}$ is defined as:
\begin{eqnarray*}
[\gamma] & \mapsto & \int_{\gamma} \omega,
\end{eqnarray*}
where $\gamma:[0,1]\to X$ is a loop represent a homotopy class of $\pi_1(X,x_0)$ with base point $x_0=\gamma(0)=\gamma(1)$.
\end{Def}


Now according to \cite{Far}, if $X$ is a CW-complex, any sigular cohomology class $\xi\in H^1(X;\mathbb{R})$ can be realised by a coninuous closed 1-form on $X$, and two closed 1-forms differ by an exact form if and only if they induce the same homomorphism of periods.




Now we have adequate volcabulary to introduce the concept of category with respect to a closed 1-form. 

\begin{Def}
Let $(X,B)$ be a finite CW pair and $\omega$ is a continuous closed 1-form on $X$, and let $N\in \mathbb{Z}$ be a positive integer and $C>0$ be a real positive constant. A subset $D\subset X$ containing $B$ is {\it $N$-movable relative to $B$ with control $C$ with respect to $\omega$} if there exists a continuous homotopy $h:D\times [0,1]\to X$ such that $h_0$ is the inclusion map, $h_t(B)\subset B$ for all $t\in [0,1]$ and for any $x\in D$, either $h_1(x)\in B$, or we have 
\begin{eqnarray*}
\int_x^{h_1(x)} \omega & \leq & -N,
\end{eqnarray*}
and for all $t\in [0,1]$, we have
\begin{eqnarray*}
\int_x^{h_t(x)} \omega & \leq & C
\end{eqnarray*}
for all $x\in D$.
\end{Def}
In this case we will simply say $D$ is {\it $(N,C)$-movable relative to $B$}. Roughly speaking a subset is $(N,C)$-movable relative to $B$ if it can be continuously deformed in the space $X$, such that any point either is pushed into $B$ or travels over distance $N$ as measured by $\omega$.

\begin{Def}
Let $(X,B)$ be a finite CW pair and $\omega$ be a continuous closed 1-form on $X$ with its cohomology class denoted as $\xi=[\omega]\in H^1(X;\mathbb{R})$. Then {\it the relative Lusternik-Schnirelmannn category with respect to $\xi$}, or $\cat(X,B,\xi)$, is defined to be the smallest integer $k$ such that there exists $C>0$ and for any integer $N>0$, there exists an open cover of $X$, $X=U\cup U_1\cup \dots \cup U_k$ such that $U_i\hookrightarrow X$ is null-homotopic in $X$ for $1\leq i\leq k$ and $U$ is $(N,C)$-movable relative to $B$.
\end{Def}

\begin{Rem}
As in the absolute case, $\cat(X,B,\xi)$ is independent of $\omega$ in the cohomology class $\xi=[\omega]$.
\end{Rem}

\begin{Rem}
When $B=\emptyset$ is empty, our $\cat (X,B,\xi)$ coincides with the controlled version of the absolute category with respect to a closed 1-form $\ccat(X,\xi)$: $\cat(X,B,\xi)=\ccat(X,\xi)$, when $B=\emptyset$. The controlled category $\ccat(X,\xi)$ was first defined in \cite{FS1}, in order to generalise the product inequality of the Lusternik-Schnirelmann category. The control is crucially used  in the proof of Theorem 1.1, however, no examples are known for which the two versions actually differ.
\end{Rem}

\begin{Rem}
When the cohomology class is trivial $\xi=0$, our category is equal to the relative version of the classical category, $\cat(X,B,\xi)=\cat(X,B)$. The notion $\cat(X,B)$ has been defined and studied in a number of papers, see for instance: \cite{Cornea}, \cite{Moyaux} and \cite{Reeken}.
\end{Rem}

This category is a homotopy invariant, the proof is analogous to the absolute case given in \cite[Section 10.2]{Far}.
\begin{Lem}
Let $\phi:(X,B)\to (X',B')$ be a relative homotopy equivalence between finite CW-complex pairs $(X,B)$ and $(X',B')$, and $\xi'\in H^{1}(X';\mathbb{R})$, $\xi=\phi^{\ast}(\xi') \in H^{1}(X;\mathbb{R})$, then

\begin{eqnarray*}
\cat(X,B,\xi) & = & \cat(X',B',\xi').
\end{eqnarray*}\eproof

\end{Lem}

\begin{center}         \section{\sc Properties of $\cat(X,B,\xi)$}                    \end{center}

We now want to prove an inequality for the relative category:

\begin{The}
Let $A\subset B\subset X$ be finite CW complexes and $\xi\in H^1(X;\mathbb{R})$ be the cohomology class of $X$, then
\begin{eqnarray*}
\cat(X,A,\xi) & \leq & \cat(X,B,\xi)+ \cat \big( B,A,i^{\ast}(\xi) \big),
\end{eqnarray*}
where the map $i^{\ast}:H^1(X;\mathbb{R})\to H^1(B;\mathbb{R})$ is induced by the inclusion map $i:B\to X$.
\end{The}

\bproof
Suppose $\cat(X,B,\xi)=k$ and $\cat(B,A,i^{\ast}(\xi))=l$, let $\omega$ be a continuous closed 1-form representing $\xi$, we need to show the existence of a real number $R>0$, such that for any $N>0$, there is an open cover of $X$ which consists of $k+l$ null-homotopic components and one $(N,R)$-movable component relative to $A$.

Firstly, we want to modify the open cover of $B$ to be open in $X$. For this we do the following trick of deformation retraction:

According to Hatcher \cite[Appendix A.2]{Hat}, there exists an open neighbourhood $N(B)$ of $B$ in $X$ such that there exists a deformation retraction $D':\overline{N(B)}\times [0,1]\to \overline{N(B)}$ rel $B$ with $D'_1(\overline{N(B)})=B$. We extend its composition with the inclusion map $\overline{N(B)}\times [0,1]\to X$ to the whole space, denoted by $D:X\times [0,1]\to X$ with $D_t|_{\overline{N(B)}}=D'_t$ for all $t$, compare to \cite[Example 0.15]{Hat}. By the compactness of $X$, there exists $K\in \mathbb{R}$ such that $\int_x^{D_1(x)}\omega< K$ for any $x\in X$.

Now according to the definition of the category, there is $C>0$ and for any integer $N$, there exist open covers $X=U\cup U_1\cup \dots \cup U_k$ and $B=V\cup V_1\cup \dots \cup V_l$, where $U_i$ and $V_j$ are null-homotopic for all $i,j$; $U$ is $(N+C+1+K, C)$-movable relative to $B$ by a homotopy $g$, and $V$ is $(N+C+2K)$-movable relative to $A$ by a homotopy $h$.

On the other hand, as $N$ varies, $\overline{N(B)}$ is not necessarily contained in $U$ for all $N>0$, therefore, let us consider the intersection $N'(B)=N(B)\cap U$ and restrict the deformation retraction to the closure of this intersection as $d=D|_{\overline{N'(B)}}: \overline{N'(B)}\times [0,1]\to X$. Note that we still have $\int_x^{d_1(x)}\omega< K$ for any $x\in \overline{N'(B)}$. Also denote by $N''(B)$ an open subset of $N'(B)$ with $N''(B)\subset \overline{N''(B)}\subset N'(B)$. In particular, $N''(B)\subset \big(d_1^{-1}(V)\cup d_1^{-1}(V_1)\cup\cdots \cup d_1^{-1}(V_l)\big)$.

Secondly, to comply with the definition of relative movability, let us modify $g:U\times [0,1]\to X$ such that points in $A$ stay in $A$ throughout the homotopy. Now according to the Lemma 3.2 below, there is an open neighbourhood $N(A)$ of $A$ in $U$ with $g_t(N(A))\subset N(B)\cap U$ for all $t\in [0,1]$. Then let $\varphi:U\to [0,1]$ be a map such that $\varphi|_A=0$ and $\varphi|_{U-N(A)}=1$. Define a continuous homotopy $g':U\times [0,1]\to X$ as
$$g'(x,t)\,\,\, =\,\,\, D(g(x,\varphi(x)t),t).$$
Then $g'_t(N(A))\subset N(B)\cap U$ as $g_t$ does for all $t$, and for any $x\in U$, either $g'_1(x)\in B$ or $\int_x^{g'_1(x)}\omega < -N-C-1$ and for all $x\in U$ and all $t\in [0,1]$, $\int_x^{g'_t(x)}\omega< C+K$.

\hfill

Now we want to show there is an open cover of $X$ modified from the ones of $X$ and $B$, namely:
$$X\,\,\, =\,\,\, (U^{\ast}\cup V^{\ast}) \cup (U_1^{\ast}\cup \dots \cup U_k^{\ast}) \cup (V_1^{\ast}\cup \dots \cup V_l^{\ast}),$$
where $U^\ast\cup V^\ast$ is $(N,R)$-movable relative to $A$ for some $R>0$ and $U_i^\ast,V_j^\ast$ are null-homotopic in $X$.

We divide the argument into three parts:

\begin{enumerate}[(i)]
 \item \textbf{Null homotopy of $V_j^{\ast}$}
To get $V_j^{\ast}$, we firstly need to modify the $V_j$'s so that they are open in $X$. Since $d$ is continuous, we have $\tilde{V}_j=d^{-1}_1(V_j)\subset N'(B)$ is open in $X$. Now we set $V_j^{\ast}=(g'_1)^{-1}(\tilde{V}_j)$ and define the null homotopy $H_j:V_j^{\ast}\times [0,1]\to X$ as 

$$ H_{j}(x,t)\,\,\, =\,\,\,  \begin{cases} g'(x,3t)              &  0           \leq t   \leq \frac{1}{3}\\
                             d(g'_1(x),3t-1)       &  \frac{1}{3} \leq t   \leq \frac{2}{3}\\
                             h_j(d_1g'_1(x), 3t-2) &  \frac{1}{2} \leq t   \leq 1,
                 \end{cases}$$
where $h^j$ is the null homotopy of $V_j$, and we see $H_j$ continuously deform $V_j^{\ast}$ to a point in $X$.

\hfill

 \item \textbf{Construction of $V^{\ast}$}
Here we want to modify $V$ and the accompanied homotopy $h$ so that the new $V^\ast$ is open in $X$ and $(N+C+K,C+1)$-movable relative to $A$ by some homotopy. Consider $V^c=B-\cup_j V_j$ in $B$, we have $d_1^{-1}(V^c)$ closed in $X$ and thus denote by $\tilde{V}^c=d_1^{-1}(V^c)\cap \overline{N''(B)}$ a closed subset in $X$. Meanwhile, let $\tilde{V}=d_1^{-1}(V)\cap N'(B) \subset N'(B)$ be open in $X$ with $\tilde{V}^c\subset \tilde{V}$. Notice that there exists a homotopy $h':\tilde{V}\times [0,1]\to X$ for $\tilde{V}$ defined as:
$$ h'(x,t)\,\,\, =\,\,\, \begin{cases} d(x,2t)        & 0            \leq t\leq \frac{1}{2}\\
                          h(d_1(x),2t-1) & \frac{1}{2}  \leq t\leq 1,
            \end{cases}$$
such that for $x\in \tilde{V}$,
$$\textrm{either } h_1'(x)\in A \textrm{ or } \int_x^{h_1'(x)}\omega < -N-C-2K+K=-N-C-K;$$
and $\int_x^{h'_t(x)} \omega < C+K$ for all $x\in \tilde{V}$ and $t\in [0,1]$.

Now according to Lemma 3.3 below, there is an open subset $V'$ of $X$ with $\tilde{V}^c \subset V'\subset \tilde{V}$ and a homotopy $H:X\times [0,1]\to X$ such that:
$$H_0(x)\,\,\, =\,\,\, x \textrm{ for all } x\in X;$$
and for all $x\in V'$, either $H_1(x)\in A$ or 
$$ \int_x^{H_1(x)} \omega \,\,\, \leq \,\,\, -N-C-K $$
and for all $x\in X$ and all $t\in [0,1]$
$$\int_x^{H_t(x)} \omega \,\,\, <\,\,\, C+1.$$
We set $V^{\ast}=(g'_1)^{-1}(V')$.

\hfill

 \item \textbf{Construction of $U^{\ast}$}
Choose slightly smaller open subsets $U_i^o\subset U_i$ such that:
$$U_i^o\subset \overline{U_i^o}\subset U_i \textrm{ and } X\subset U\cup U_1^o\cup \dots \cup U_k^o,$$
then we define
$$U^{\ast}\,\,\, =\,\,\, X-\left(\big(\bigcup_{i=1}^k \overline{U^o}_i\big) \cup (g'_1)^{-1}(\overline{N''(B)})\right).$$

Define the homotopy $G:\big( U^{\ast}\cup V^{\ast}\big)\times [0,1] \to X$ as:
\begin{displaymath}
G(x,t)\,\,\, =\,\,\,      \begin{cases}
                g'(x,2t)     &  0           \leq t  \leq \frac{1}{2}\\
                H(g'_1(x),2t-1)    &  \frac{1}{2} \leq t  \leq 1  
                \end{cases}.
\end{displaymath}
It is easy to see that $G_t(A)\subset A$ for all $t\in [0,1]$ as both $g'$ and $H$ are built with this feature. For $x\in U^{\ast}$ it will travel over distance $N$ as:
$$\int_x^{G_1(x)}\omega=\int_x^{g'_1(x)} \omega + \int_{g'_1(x)}^{H_1(g'_1(x))}\omega\leq (-N-C-1)+(C+1) = -N.$$

Similarly, for $x\in V^\ast=(g'_1)^{-1}(V')$, after discounting the effect of $g'$ and returning into $V'\subseteq N(B)$, $H$ either pushes the point into $A$ or travel over distance $N$ as
$$\int_x^{G_1(x)}\omega=\int_x^{g'_1(x)} \omega + \int_{g'_1(x)}^{H_1(g'_1(x))}\omega\leq C+K+(-N-C-K) = -N.$$
Also for all $t\in [0,1]$ and $x\in U^\ast\cup V^\ast$, $\int_x^{G_t(x)} \omega < 2C+2K+1$.
\end{enumerate}

Finally, let us set $U_i^{\ast}=U_i$ unchanged, then $X$ is covered as:
$$X\,\,\, =\,\,\, (U^{\ast}\cup V^{\ast}) \cup (U_1^{\ast}\cup \dots \cup U_k^{\ast}) \cup (V_1^{\ast}\cup \dots \cup V_l^{\ast}).$$
This is true as $(g'_1)^{-1}(\overline{N''(B)})$ is covered by $V^\ast$ and $V_j^\ast$:
$$(g'_1)^{-1}(\overline{N''(B)})\,\,\, \subset\,\,\, V^\ast\cup V_1^\ast\cup\dots \cup V_l^\ast,$$
where
\begin{eqnarray*}
\overline{N''(B)} & \subset & d_{\frac{1}{2}}^{-1}(V^c)\cup d_1^{-1}(V_1)\cup \dots \cup d_1^{-1}(V_l)\\
                  & \subset & V'\cup \tilde{V}_1\cup \dots \cup \tilde{V}_l;
\end{eqnarray*}
and $\{U_i^\ast\}$ covers the rest of $X$.

Now $U^\ast\cup V^\ast$ is $(N,2C+2K+1)$-movable relative to $A$ and the other components are all null-homotopic.
\eproof

\begin{Lem}
With the notations as in the proof of Theorem 3.1, there exists open neighbourhood $N(A)$ of $A$ in $X$ with $N(A)\subset N(B)\cap U$, such that $g_t(N(A))\subset N(B)\cap U$ for all $0\leq t\leq 1$.
\end{Lem}

\bproof
Given $g: U\times [0,1]\to X$, we have $g_t(a)\in B\subset N(B)\cap U$ for any $a\in A$, according to the hypothesis. For such point $(a,t)\in A\times [0,1]$, by the continuity of $g$, we can find some neighbourhood $N^t(a)\times (t-\delta_t, t+\delta_t)$ of $(a,t)$ in $X$ for small $\delta_t$, such that $g(a',t')\in N_{\epsilon}(B)\cup U$ for all $(a',t')\in N^t(a)\times (t-\delta_t, t+\delta_t)$.

Now because of the compactness of $[0,1]$, there exists $t_1,\dots,t_n$ such that $[0,1]=\bigcup_{i=1}^n(t_i-\delta_i,t_i+\delta_i)$. Set
$$N(a)=\bigcap_{i=1}^n N^{t_i}(a),$$
we claim $g(N(a)\times [0,1])\in N(B)\cap U$.


Now define 
$$N(A)=\bigcup_{a\in A} N(a).$$

We can see $N(A)\subset N(B)\cap U$ and $g_t(N_\delta(A))\subset N(B)\cap U$ for all $t\in [0,1]$.
\eproof

The following lemma is a convenient generalisation of Lemma 10.1 in \cite{FS1}, stating that the homotopy for a movable subset can be extended to the whole space $X$ with the control $C+1$, and the proof follows essentially the same argument as in \cite{FS1}.

\begin{Lem}
Let $\omega$ be a continuous closed 1-form on a finite CW complex $X$. Let $B\subset X$ be a subcomplex. Suppose further that there exists a $C>0$ and for any integer $N>0$, we have an open subset $U$ of $X$ containing $B$ and $U$ is $(N,C)$-movable with respect to $B$. Then for any given closed subset $W\subset U$ with $B\subset W$, there exists an open set $U'$ with $W\subset U'\subset U$ and a homotopy $H:X\times [0,1]\to X$ satisfying the following:

\begin{enumerate}
 \item $H_0(x)=x$ for all $x\in X$ and $H_t(B)\subset B$ for all $t\in [0,1]$; 
 \item For any $x\in U'$ one has either $H_1(x)\in B$ or $\int_x^{H_1(x)} \omega < -N$;
 \item For any $x\in X$ and $t\in [0,1]$, $\int_x^{H_t(x)} \omega < C+1$.
\end{enumerate}\eproof
\end{Lem}

If $A=\emptyset$ is empty, we get the following corollary:
\begin{Cor}
Let $(X,B)$ be a finite CW pair and $\xi\in H^1(X;\mathbb{R})$, then
$$\cat(X,\xi)\,\,\, \leq\,\,\, \cat(X,B,\xi)+ \cat \big( B,i^{\ast}(\xi) \big).$$\eproof
\end{Cor}

We can also derive a similar inequality for the category of a product of CW-complex pairs, compare with \cite{FS1}.

\begin{The}
Let $(X,B),(Y,D)$ be two CW pairs, $\xi_{X}\in H^{1}(X;\mathbb{R})$ and $\xi_{Y}\in H^{1}(Y;\mathbb{R})$ be the cohomology classes on $X$ and $Y$, respectively. Suppose also

\begin{displaymath}
\cat(X,B,\xi_{X})\,\,\, >\,\,\, 0 \qquad \textrm{or} \qquad \cat(Y,D,\xi_{Y})\,\,\, >\,\,\, 0,
\end{displaymath}

Then 
\begin{displaymath}
\cat((X,B)\times (Y,D),\xi)\,\,\, \leq\,\,\, \cat(X,B,\xi_{X})+\cat(Y,D,\xi_Y)-1,
\end{displaymath}
with $\xi=\xi_X \times 1+1\times \xi_Y$.\eproof
\end{The}

We now want to provide a cohomology lower bound for $\cat(X,B,\xi)$ similar to the one in \cite{FS1}. Let us begin with some basic notions.

For a CW complex $X$ and a continuous closed 1-form $\omega$, we have a regular covering space $p:\tilde{X}\to X$ correspond to the kernel of the cohomology class $\xi=[\omega]\in H^1(X;\mathbb{R})$. The covering transformation group is $H\simeq \mathbb{Z}^r=\pi_1(X)/\ker(\xi)$. Then the cohomology class of the pullback of $\omega$ is trivial in the covering, $[p^{\ast}\omega]=0\in H^1(\tilde{X};\mathbb{C})$, that is, there exists a real function $f:\tilde{X}\to \mathbb{R}$ such that $\ud f=p^\ast \omega$.

\begin{Def}
A subset $O\subset X$ is called a {\it neighbourhood of infinity} in $\tilde{X}$ with respect to a cohomology class $\xi\in H^1(X;\mathbb{R})$, if $O$ contains the set $\{x\in \tilde{X}: f(x)<c\}$ for a real number $c\in \mathbb{R}$. Here $f:\tilde{X}\to \mathbb{R}$ is a real function obtained by pulling back a closed 1-form $\omega$ with $[\omega]=\xi\in H^1(X;\mathbb{C})$.
\end{Def}

Notice the definition of a neighbourhood of infinity $O$ is independent of the choice of real functions. For a more detailed exposition of this concept, we refer to \cite{FS3}, in particular Lemma 3 in Section 3.


\begin{Def}
Let $(X,B)$ be a finite CW complex pair and $\omega$ be a continuous closed 1-form on $X$. Suppose $p:\tilde{X}\to X$ is a regular covering corresponding to $\ker(\xi)$ where $\xi=[\omega]\in H^1(X)$ is the cohomology class of $\omega$. Then a homology class $z\in H_i(\tilde{X},\tilde{B})$ is {\it movable to infinity with respect to $\xi$}, if in any neighborhood $O$ of infinity with respect to $\xi$, there exists a relative homology class in $H_i(O,O\cap \tilde{B})$ whose image is $z$ under the map  $H_i(O,O\cap \tilde{B})\to H_i(\tilde{X},\tilde{B})$ induced by inclusion.
\end{Def}



\begin{Not}
Let $H=H_1(X;\mathbb{Z})/\ker(\xi)$, denote $\mathcal{V}_\xi=(\mathbb{C}^\ast)^r=\homo(H,\mathbb{C}^\ast)$, which we can think of as the variety of all complex flat line bundles $L$ over $X$ such that the induced flat line bundle $p^\ast L$ on $\tilde{X}$ is trivial.
\end{Not}


\begin{Def}
In $\mathcal{V}_\xi$ a bundle $L$ is called {\it $\xi$-transcendental} if the monodromy $\mon_L:\mathbb{Z}[H]\to \mathbb{C}$ is injective, and {\it $\xi$-algebraic} if not.
\end{Def}

The following two assertions are the relative versions of Proposition 6.5 and Theorem 4 in \cite{FS1}, their validity follows similar algebraic arguments provided in \cite{FS1}.

\begin{Pro}
Suppose $L\in \mathcal{V}_\xi$ is $\xi$-transcendental, and $v\in H^q(X,B;L)$ is a non-zero cohomology class. Then there exists a homology class $z\in H_q(\tilde{X},\tilde{B};\mathbb{C})$ with $v\frown p_\ast(z)\neq 0$.
\eproof
\end{Pro}

\begin{The}
Suppose a flat line bundle is $\xi$-transcendental and there is cohomology class $v\in H^q(X,B;L)$ with $v\frown p_\ast(z)\neq 0$ for some $z\in H_q(\tilde{X},\tilde{B};\mathbb{C})$ and $p_\ast(z)\in H_q(X,B;L^\ast)$, where $L^\ast$ is the dual bundle of $L$. Then $z$ is not movable to infinity with respect to $\xi$.
\eproof
\end{The}

We now state the cohomology estimate of the category:

\begin{The}
Suppose $L\in \mathcal{V}_\xi$ is $\xi$-transcendental, $v_0\in H^{d_{0}}(X,B;L)$ and $v_i\in H^{d_i}(X;\mathbb{C})$ for $i=1,\dots,k$ such that $d_i>0$ and
\begin{equation}
\label{eq:cup}
v_0\smile v_1\smile \cdots \smile v_k \,\,\, \ne \,\,\, 0\in H^{d}(X,B;L),
\end{equation}
with $d=\Sigma_{i} d_i$, then
\begin{displaymath}
\cat(X,B,\xi)\,\,\, >\,\,\, k.
\end{displaymath}
\end{The}

The maximal such $k$ gives a lower bound for $\cat(X,B,\xi)$, and it gives a cup length estimate for $\cat(X,B,\xi)$.

\hfill

\bproof
Let $v=v_1\smile \cdots \smile v_k$, according to (\ref{eq:cup}) and Proposition 3.10, we can find a homology class $z\in H_d(\tilde{X},\tilde{B};\mathbb{C})$ such that
$$(v_0\smile v)\frown p_\ast(z)\,\,\, \neq\,\,\, 0.$$

Fix such a homology class $z\in H_d(\tilde{X},\tilde{B};\mathbb{C})$, then it is possible to choose a compact polyhedron $K\subset \tilde{X}$ such that $z$ is the image of some homology class in $H_d(K,\tilde{B}\cap K;\mathbb{C})$ under the inclusion-induced map $i_\ast: H_d(K,\tilde{B}\cap K;\mathbb{C})\to H_d(\tilde{X},\tilde{B};\mathbb{C})$. We denote this homology class $z'\in H_d(K,\tilde{B}\cap K;\mathbb{C})$. Now we assert the existence of a neighbourhood of infinity $O_\infty$ which possesses the following property: if the image of a homology class under the map $H_\ast(K,\tilde{B}\cap K;\mathbb{C})\to H_\ast(\tilde{X},\tilde{B};\mathbb{C})$ has a preimage in $H_\ast(O_\infty,O_\infty\cap \tilde{B};\mathbb{C})$, then it is movable to infinity. Indeed, let $O=f^{-1}((-\infty,0])\subset \tilde{X}$ be a neighbourhood of infinity, and $g:\tilde{X}\to \tilde{X}$ be a covering transformation such that $\xi(g)<0$. Then

$$V_g\,\,\, =\,\,\, \im [H_\ast(gO,gO\cap \tilde{B};\mathbb{C})\to H_\ast(\tilde{X},\tilde{B};\mathbb{C})]\cap \im [H_\ast(K,\tilde{B}\cap K;\mathbb{C})\to H_\ast(\tilde{X},\tilde{B};\mathbb{C})]$$
is a finite dimensional complex vector space.


We get a chain of finite dimensional vector spaces:
$$\cdots \subset V_{g^n}\subset \cdots \subset V_{g^2}\subset V_g\subset V$$
which stabilises after finitely many terms. Subsequently, there exists a sufficiently large $N>0$ such that $V_{g^n}= V_{g^N}$ for any $n\geq N$. Therefore, fix such a $N$ and the subset $O_\infty=g^NO$ will work.

So let us have such a neighbourhood $O_\infty$, then the pullback function $f:\tilde{X}\to \mathbb{R}$ of $\omega$ with $p^\ast\omega=\ud f$ gives values to points in $K$ and $O_\infty$. In particular, we have $f(K)\subset [a,b]$ and $O_\infty\supset f^{-1}(-\infty,c)$, for some $c<a<b$. Note that $c<a$ is always possible by increasing $N$ if necessary.

Now assume the statement is false, then $\cat (X,B,\xi)\leq k$, in particular, for $N> b-c$ and some $C>0$, there exists an open cover of $X$:
$$X\,\,\, =\,\,\, U\cup U_1\cup \cdots \cup U_k,$$
where $U_i \hookrightarrow X$ is null-homotopic and $U$ is $(N,C)$-movable relative to $B$.

Now observe that $v_i\in H^{d_i}(X;\mathbb{C})$ can be pulled back to some $u_i\in H^{d_i}(X,U_i;\mathbb{C})$ because of the null-homotopy of $U_i$.

Therefore, by naturality of the cup product, $v=j^\ast(u)\in H^{d-d_0}(X;\mathbb{C})$ for some $u=u_1\smile \cdots \smile u_k\in H^{d-d_0}(X,U_1\cup \cdots \cup U_k;\mathbb{C})$, where $j^\ast$ is induced by inclusion $j:(X,\emptyset)\to (X,U_1\cup \cdots \cup U_k)$.

Let $w$ be the image of $p^\ast(u)$ via the inclusion-induced map
$$i_1^\ast: H^{d-d_0}(\tilde{X},\tilde{U_1}\cup \cdots \cup \tilde{U_k};\mathbb{C})\to H^{d-d_0}(K,(\tilde{U_1}\cup \cdots \cup \tilde{U_k})\cap K;\mathbb{C}),$$ and restrict the lift $(\tilde{X},\emptyset)\to (\tilde{X},\tilde{U}_1\cup\cdots \cup\tilde{U}_k)$ of $j$ to $K$ as: $$\tilde{\jmath}:(K,\emptyset)\to (K,(\tilde{U}_1\cup\cdots \cup\tilde{U}_k)\cap K),$$
then
$$\tilde{\jmath}^\ast w\frown z' \,\,\, \in\,\,\, H_{d_0}(K,\tilde{B}\cap K),$$
where $\tilde{\jmath}^\ast w\in H^{d-d_0}(K; \mathbb{C})$ and $z'\in H_d(K,\tilde{B}\cap K;\mathbb{C})$. Notice $\tilde{\jmath}^\ast w\frown z'\ne 0$ as by naturality of the cap product, see \cite[Lemma 5.6.16, pp~254]{Spanier},
\begin{eqnarray}
i_\ast (\tilde{\jmath}^\ast w\frown z') 
& = & i_\ast (\tilde{\jmath}^\ast i_1^\ast(p^\ast(u))\frown z') \,\,\,=\,\,\, i_\ast ( (i_1\tilde{\jmath})^\ast(p^\ast(u))\frown z')\nonumber\\
& = & p^\ast(u)\frown i_{2\ast}(z')                             \,\,\,=\,\,\, p^\ast(u)\frown j_{1\ast}i_\ast(z')\nonumber\\
& = & j_2^\ast p^\ast(u)\frown i_\ast(z')                       \,\,\,=\,\,\, (p j_2)^\ast (u) \frown z\nonumber\\
& = & (j p)^\ast (u) \frown z                                   \,\,\,=\,\,\, p^\ast (j^\ast(u))\frown z \nonumber\\
& = & p^\ast (v) \frown z \nonumber,
\end{eqnarray}
which is non-trivial according to our hypothesis. Here $j_2^\ast: H_d(\tilde{X},\tilde{U}_1\cup\dots\cup\tilde{U}_k)\to H_d(\tilde{X})$ and $i_{2\ast}:H_d(K,\tilde{B}\cap K;\mathbb{C})\to H_d(\tilde{X},\tilde{B}\cup \big(\tilde{U}_1\cup\dots\cup\tilde{U}_k\big))$ are induced by the inclusion map $j_2$ and $j_1 i$ with $j_1:(\tilde{X},\tilde{B})\to (\tilde{X},\tilde{B}\cup \big(\tilde{U}_1\cup\dots\cup\tilde{U}_k\big))$, respectively.


Again by naturality of the cap product,
$$i'_\ast(\tilde{\jmath}^\ast w\frown z')\,\,\, =\,\,\, w\frown \bar{\imath}_\ast(z),$$

where $i'_\ast$ is induced by $i':(K,\tilde{B}\cap K)\to (K,\tilde{U}\cap K)$ and $\bar{\imath}_\ast$ is from 
$$\bar{\imath}:(K,\tilde{B}\cap K)\to (K,\left((\tilde{U}_1\cup\cdots \cup\tilde{U}_k)\cap K\right)\cup \left(\tilde{U}\cap K\right))=(K,K).$$

Therefore, $\bar{\imath}_\ast(z)=0\in H_{d_0}(K,K)=0$, and $i'_\ast(\tilde{\jmath}^\ast w\frown z)\in H_{d_0}(K,\tilde{U}\cap K)$ is trivial. Consequently, the exact sequence
$$\cdots \to H_{d_0}(\tilde{U}\cap K, \tilde{B}\cap K)\to H_{d_0}(K,\tilde{B}\cap K)\stackrel{i'_\ast}{\to} H_{d_0}(K,\tilde{U}\cap K)\to \cdots$$
indicates the existence of a nontrivial preimage $z_0$ of $\tilde{\jmath}^\ast w\frown z$ in $H_{d_0}(\tilde{U}\cap K, \tilde{B}\cap K)$.

Now for the $(N,C)$-movable open subset $U$ in $X$, its lift $\tilde{U}$ in $\tilde{X}$ has a homotopy $h:(\tilde{U},\tilde{B})\times [0,1]\to (\tilde{X},\tilde{B})$ starting with the inclusion, and $fh_1(x)-f(x)\leq -N $, hence $(h_0)_\ast z_0 = (h_1)_\ast z_0\in H_{d_0}(\tilde{X},\tilde{B})$. But since $h_1:(\tilde{U}\cap K,\tilde{B}\cap K)$ factors through $(\tilde{U}\cap K,\tilde{B}\cap K)\to (O_\infty,O_\infty\cap \tilde{B})$ due to our choice of $N>b-c$ for $\tilde{U}$, there exists a homology class in $H_{d_0}(O_\infty,O_\infty\cap \tilde{B})$ that maps to $(h_1)_\ast(z_0)$. In other words,
\begin{eqnarray}
v_0\frown p_\ast(h_\ast z_0) & = & v_0\frown p_\ast (i_\ast(\bar{\jmath}^\ast w\frown z'))\,\,\, =\,\,\, v_0\frown p_\ast(p^\ast(v)\frown z)\nonumber\\
                                & = & (v_0\smile v)\frown p_\ast(z)\,\,\, \neq\,\,\, 0,\nonumber
\end{eqnarray}
in contradiction to Theorem 3.11.
\eproof

\begin{center}                \section{\sc Homoclinic cycles and critical points}                                               \end{center}
In this section, let $M$ be a smooth compact manifold with boundary $\partial M$. We relate the invariant $\cat (M,B,\xi)$ to critical points of smooth closed 1-form $\omega$ representing a cohomology class $\xi\in H^1(X;\mathbb{R})$. Here $B\subset X$ and $\omega$ are related as follows:

Let $\rho:\bar{M}\to M$ be a regular covering space of $M$ with $\pi_1(\bar{M})=\ker([\omega])$, then there exists a smooth function $f:\bar{M}\to \mathbb{R}$ with $\ud f=\rho^\ast(\omega)$. We equip the boundary $\partial M$ of the $M$ with a tubular neighbourhood structure $\partial M\times [0,1)\subset M$ and lift it up to a tubular neighbourhood $\partial \bar{M}\times [0,1)\subset \bar{M}$ in the covering. Points in the tubular neighbourhood will be denoted as $(x,t)\in \partial \bar{M}\times [0,1)$. Fixing this neighbourhood, for $x\in \partial \bar{M}$ we get a well-defined partial derivative $\dfrac{\partial f}{\partial t}|_{(x,0)}\in \mathbb{R}$. Notice it is equivariant with respect to the action of the transformation group of $\rho$, which implies a smooth map $\dfrac{\partial f}{\partial t}: \partial M\to \mathbb{R}$ on boundary of the base manifold.
\begin{Def}
Given a fixed inner collaring $\partial M\times [0,1)$ of the boundary $\partial M$ of $M$, {\it the exit set $B$ of $\omega$} is defined as:
$$B \,\,\, = \,\,\ \{x\in \partial M: -\dfrac{\partial f}{\partial t}|_{(x,0)}\leq 0 \}.$$
\end{Def}

By a gradient of $\omega$ we mean a vector field $v$ which is dual to $\omega$ with respect to some Riemannian metric.

\begin{Not}
Let $\Phi: \Delta \to M$ be the negative gradient flow of a gradient vector field $v$ of $\omega$, where $\Delta\subset M\times \mathbb{R}$.
\end{Not}

We want to have that $B$ is the set where the negative flow 'exits' the manifold. For this we need some restriction on $\omega$ and the gradients.

We describe the conditions in terms of the pullback $\ud f=\rho^\ast \omega$:
\begin{description}
\item[Assumptions on $\omega$ on $\partial M$]
\end{description}

\begin{enumerate}[]
\item[$\mathbf{A1}$] The function $f$ has no critical point on $\partial \bar{M}$. Without loss of generality we assume that $f$ has no critical points in the entire collaring $\partial \bar{M}\times [0,1)$.

\item[$\mathbf{A2}$] The partial derivative $\frac{\partial f}{\partial t}$, where $t$ is the coordinate for $[0,1)$, is a smooth function on $\partial \bar{M}\times \{0\}$ and hence on $\partial M$, and zero is a regular value of $\frac{\partial f}{\partial t}(x,0)$. Denote by $\Gamma=\{ x\in \partial M: \frac{\partial f}{\partial t}(x,0)=0\}$, this is equivalent to say $\Gamma$ is a $1$-codimensional closed submanifold of $\partial M$.

\item[$\mathbf{A3}$] Fix a tubular collaring of $\Gamma$ in $\partial M$, $\Gamma\times [-1,1]\subset \partial M$, with $\Gamma\times [-1,0]\subset B$. So if a point lies in the cubical neighbourhood of $\Gamma$ in $M$, we write it in local coordinates:
$$(x,s,t)\,\,\, \in\,\,\, \Gamma\times [-1,1]\times [0,1),$$
where $x=(x_1,\cdots, x_{m-2})$, then we assume
$$\frac{\partial f}{\partial s}(x,0,0)\,\,\, >\,\,\, 0.$$

\end{enumerate}

Notice that the conditions $\mathbf{A1},\mathbf{A2}$ and $\mathbf{A3}$ do not depend on the particular choice of collarings. However, in order to get that $B$ serves as the exit set for the negative gradient flow, we need a restriction on the gradients. We formalise the idea by the following notion:

\begin{Def}
Let $\omega$ be a closed 1-form that satisfies $\mathbf{A1},\mathbf{A2}$ and $\mathbf{A3}$. A gradient $v$ of $\omega$ is called {\it transverse on $(\partial M,B)$} if the Riemannian metric is the product metric on $\Gamma$ and on $\partial M$ with respect to the same tubular neighbourhoods as in $\mathbf{A2}$ and $\mathbf{A3}$.
\end{Def}

With such a restriction on gradients, we get the following lemma on the 'timing' of the moment at which each point reaches $B$:
\begin{Lem}
Let $v$ be a gradient of $\omega$ transverse on $(\partial M,B)$ and denote by
$$U_B=\{x\in M: \textrm{ there exists } t\in \mathbb{R}, \textrm{ such that } x\cdot t\in B\},$$
where $x\cdot t$ is a shorthand notation of the negative gradient flow $\Phi(x,t)$ for each $x$ and $t$. Then the function $\beta: U_B\to \mathbb{R}$ defined as $\beta(x)=\min \{t:x\cdot t\in B\}$ is continuous, and $U_B$ is open in $M$.\eproof
\end{Lem}
The proof is given in the first author's thesis \cite[Section 1.2]{I}. The crucial point is that interior points of $M$ cannot reach $\Gamma$ under the flow if $v$ is transverse on $(\partial M,B)$.

Now let us recall the definition of homoclinic cycle which is a generalisation of homoclinic orbit. Here we implicitly assume that $\omega$ has only finitely many crticial points. For a more general treatment see \cite{Latschev}.

\begin{Def}
A sequence of trajectories $\{\gamma_i(t):\mathbb{R}\to M \}_{1\leq i\leq n}$ on a manifold $M$ is called a {\it homoclinic cycle of length $n$} if for each $\gamma_i$ its limit $\lim_{t\to \pm \infty} \gamma_i(t)$ exists and the following is satisfied:
$$ \lim_{t\to +\infty} \gamma_i(t)=\lim_{t\to -\infty} \gamma_{i+1}(t) \textrm{ for } 1\leq i \leq n-1, \textrm{ and } \lim_{t\to +\infty} \gamma_n(t)=\lim_{t\to -\infty} \gamma_1(t).$$
\end{Def}

\begin{Def}
A trajectory $\gamma$ is said to have {\it displacement $N$ by $\omega$} if its integral with respect to $\omega$ equals $N$:
\begin{eqnarray*}
\int_{\gamma} \omega & = & N,
\end{eqnarray*}
a homoclinic cycle $\{ \gamma_i \}$ has {\it displacement $N$ by $\omega$} if
\begin{eqnarray*}
\sum_i \int_{\gamma_i} \omega & = & N.
\end{eqnarray*}
\end{Def}

\begin{The}
Let $M$ be a smooth compact  manifold with boundary $\partial M$, and $\omega$ be a closed 1-form on $M$ with exit set $B\subset \partial M$ satisfying assumptions $\mathbf{A1,A2}$ and $\mathbf{A3}$ below. If the number of critical points of $\omega$ is less than $\cat(M,B,[\omega])$, then any gradient of $\omega$ transverse on $\partial M$ contains at least one homoclinic cycle.
\end{The}




\bproof
For any real number $N>0$, assume there is a gradient of $\omega$ transverse on $(\partial M,B)$ without homoclinic cycle of displacement less than $N$. For some $C>0$ and any such $N>0$ we need to show the existence of an open cover $M=U\cup U_1\cup \dots \cup U_k$ according to the definition of $\cat(M,B,\xi)$, where $\xi=[\omega]\in H^1(M;\mathbb{R})$ is the cohomology class of $\omega$.

The idea is to use the negative gradient flow as the prototype for the homotopies and partition the manifold according to the destination of each point travelling along its flow line.

Because the homotopy is modified from the negative gradient flow, the integral $\int \omega\leq 0$ is always non-positive along the trajectories, so we can choose $C=0$. Let us fix $N>0$, we want to construct an open cover of $M$ as
$$M=U\cup U_1\cup \dots \cup U_k.$$

We firstly define $U$ as the open subset of all the points either reach $B$ in finite time or travel over displacement $N$ in the negative direction:
$$U=\{x\in M: \textrm{ there exists some } t_x>0 \textrm{ such that either } x\cdot t_x\in B, \textrm{ or } \int_x^{x\cdot t_x} \omega < -N\}.$$

Secondly, for $U_i$, we first need a so-called {\it gradient-convex neighbourhood} $V_i$ for each critical points $p_i$, in order to construct open subsets. For each critical point $p_i$, the gradient-convex neighbourhood $V_i$ is a small closed disc containing $p_i$, such that the points on the boundary of $V_i$ who are leaving $V_i$ under the negative gradient flow have to travel over displacement $N$ before returning to $\intl V_i$. The existence of $V_i$ is derived from the no homoclinic cycle condition in the hypothesis, for a detailed argument see \cite{Far1} and \cite{Latschev}. Then we define $U_i$ for each $p_i$ as follows:
$$U_i=\{x\in M: x\cdot t_x \in \intl V_i \textrm{ for some } t_x\in \mathbb{R} \textrm{ and } \int_x^{x\cdot t_x} \omega >-N\}.$$
The null homotopy of $U_i$ can also be found proved in \cite{Far1} and \cite{Latschev}.

Now we are left to show the movability of $U$. The subset $U$ is open since it is the union of two open subsets, namely $\{x\in M: \int_x^{x\cdot t_x} \omega < -N \textrm{ for some } t_x>0, \textrm{ where } N>0\}$ and $\{x\in M: x\cdot t_x\in B \textrm{ for some } t_x>0\}$, they are both open by Lemma 4.4 and continuity of the flow.



According to the construction, for each $x\in U$, there exists $t_x\in \mathbb{R}$, such that either $x\cdot t_x\in B$ or $\int_x^{x\cdot t_x} \omega =-N$, by Lemma 4.4 and the Implicit Function Theorem. Moreover, the map $x\to t_x$ is a real continuous function on $U$. Therefore, we can define the homotopy $h:U\times [0,1]\to M$ as
$$h(x,\tau)=x\cdot (\tau t_x).$$
This proves Theorem 4.7, hence Theorem 1.2.
\eproof

Notice that the the homotopy $h:U\times [0,1]\to M$ in the proof above fiexes $B$, so we could consider modifying the definition of $\cat(X,B,\xi)$ by demanding the homotopy to fix $B$. However, this leads to the same number, see \cite{I}.

\end{document}